\documentclass[10pt]{amsart}

\usepackage[all]{xy}
\usepackage{latexsym}
\usepackage{amsmath, mathtools}
\usepackage{amsthm}
\usepackage{amssymb}
\usepackage{a4wide}
\usepackage[utf8]{inputenc}
\usepackage{dsfont}
\usepackage{multirow}
\usepackage{slashbox}
\usepackage{pdfsync}
\usepackage{enumitem}
\usepackage{tikz-cd}
\usepackage{hyperref}

\swapnumbers
\newtheorem{thm}[subsubsection]{Theorem}
\newtheorem{dei}[subsubsection]{Definition}

\newtheorem{lem}[subsubsection]{Lemma}

\theoremstyle{definition}
\newtheorem*{rem}{Remark}

\newtheorem*{pf}{\sc Proof}

\newtheorem*{ex}{Examples}

\newcommand{\B}{\mathrm{B}}
\newcommand{\C}{\mathcal{C}}

\newcommand{\F}{\mathcal{F}}

\newcommand{\Po}{\mathcal{P}}
\newcommand{\Poa}{\Po^{\ash}}
\newcommand{\So}{\mathbb{S}}

\newcommand{\cqfd}{\ \hfill \square}

\newcommand{\ash}{\textrm{!`}}

\newcommand{\epi}{\twoheadrightarrow}
\newcommand{\mono}{\rightarrowtail}

\newcommand{\uAs}{u\mathcal{A}s}

\newcommand{\id}{\operatorname{id}}

\newcommand{\Hom}{\operatorname{Hom}}

\newcommand{\un}{\vcenter{
\xymatrix@M=0pt@R=3pt@C=3pt{&&\\
& \ar@{{*}}[u] \ar@{-}[d]\\
& &}}}

\newcommand{\ass}{\vcenter{ \xymatrix@M=0pt@R=4pt@C=4pt{
      & & & & \\ \ar@{-}[ddrr] & & \ar@{-}[dl] & & \ar@{-}[ddll]\\
      &  & &  & \\
      & & \ar@{-}[d] & &\\
      & & }}
- 
\vcenter{
    \xymatrix@M=0pt@R=4pt@C=4pt{& & & & \\
      \ar@{-}[drdr] & &\ar@{-}[dr] & & \ar@{-}[dldl]  \\
      & & & & \\
      & &\ar@{-}[d] & & \\
      & & }}}

\newcommand{\mun}{\vcenter{ \xymatrix@M=0pt@R=4pt@C=4pt{
      & & & & \\
      & & & &\\
      & \ar@{{*}}[u] \ar@{-}[dr] & & \ar@{-}[dl] & \\
      & & \ar@{-}[d] & &\\
      & & }}}

\newcommand{\munsplit}{\vcenter{ \xymatrix@M=0pt@R=4pt@C=4pt{
      & & & & \\
      &\ar@{{*}}[u] & & &\\
      & \ar@{-}[dr] & & \ar@{-}[dl] & \\
      & & *{} \ar@{-}[d] & &\\
      & & }}}

\newcommand{\lmun}{\vcenter{\xymatrix@M=0pt@R=7pt@C=10pt{
      & & & & \\
      & & & &\\
      & \ar@{{*}}[u] \ar@{-}[dr] & & \ar@{-}[dl] & \\
      & & \ar@{-}[d] & &\\
      & & }}}

\newcommand{\unm}{\vcenter{
    \xymatrix@M=0pt@R=4pt@C=4pt{& & & & \\
      & & & &\\
      & \ar@{-}[dr] & & \ar@{{*}}[u] \ar@{-}[dl] & \\
      & &\ar@{-}[d] & & \\
      & & }}}

\newcommand{\unmsplit}{\vcenter{
    \xymatrix@M=0pt@R=4pt@C=4pt{& & & & \\
      & & & \ar@{{*}}[u] &\\
      & \ar@{-}[dr] & & \ar@{-}[dl] & \\
      & &\ar@{-}[d] & & \\
      & & }}}

\newcommand{\mtwo}{\vcenter{\xymatrix@M=0pt@R=7pt@C=10pt{
    &\ar@{-}[dr] & & \ar@{-}[dl] \\
    & &\ar@{-}[d] & \\
    & & & }}}

\newcommand{\smtwo}{\vcenter{\xymatrix@M=0pt@R=4pt@C=4pt{
    &\ar@{-}[dr] & & \ar@{-}[dl] \\
    & &\ar@{-}[d] & \\
    & & & }}}

\newcommand{\mthree}{\vcenter{
    \xymatrix@M=0pt@R=4pt@C=4pt{
      \ar@{-}[dr] & \ar@{-}[d]& \ar@{-}[dl] &\\
      &\ar@{-}[dd] &  &\\
      & & & \\
      & & & }}}
\newcommand{\mthreec}{\vcenter{
    \xymatrix{
       *{}\ar@{-}[rd] &*{} \ar@{-}[d]& *{}\ar@{-}[ld] & \\
       &*{}\ar@{-}[d] & &  \\
       & & &
    }}}
\newcommand{\mfivectwo}{\vcenter{ \xymatrix@M=0pt@R=7pt@C=10pt{
    & & & \\
    & & & & \\
    \ar@{{*}}[u] \ar@{-}[drr] & \ar@{-}[dr] & \ar@{-}[d] &
    \ar@{-}[dl] \ar@{{*}}[u]  & \ar@{-}[dll] \\
    & & \ar@{-}[d] \\
    & &
    }}}

\newcommand{\smfivectwo}{\vcenter{\xymatrix@M=0pt@R=4pt@C=4pt{
    & & & \\
    & & & & \\
    \ar@{{*}}[u] \ar@{-}[drr] & \ar@{-}[dr] & \ar@{-}[d] &
    \ar@{-}[dl] \ar@{{*}}[u]  & \ar@{-}[dll] \\
    & & \ar@{-}[d] \\
    & &
    }}}

\newcommand{\ident}{\vcenter{\xymatrix@M=0pt@R=7pt@C=10pt{
  \ar@{-}[dd]& \\
  & \\
  &
}}}

\newcommand{\draftnote}[1]{}

\author{Joseph Hirsh and Joan Mill\`es}
\thanks{The first author was supported by a
  National Science Foundation Graduate Research Fellowship.}
\thanks{The second was supported by the ANR grant JCJC06 OBTH}

\title{Corrigendum for the article\\
``Curved Koszul duality theory''}


\begin{document}

\begin{abstract}
In this corrigendum, we explain and correct a mistake in our article ``Curved Koszul duality theory'' \cite{HirshMilles}. Our definitions of morphisms between semi-augmented properads and between curved coproperads have to be modified.
\end{abstract}

\maketitle

\tableofcontents

\section*{Introduction}
In the \cite{HirshMilles}, we define the notion of \emph{semi-augmented} dg properad and we write that a morphism of such dg properads is a morphism $f$ of dg properads which commutes the semi-augmentations $\varepsilon$ and $\varepsilon'$ that is $\varepsilon' \cdot f = \varepsilon$. 

We should not require this last condition. Indeed we provide in the original article two maps of semi-coaugmented dg properads which are the bar-cobar resolution $g_\pi : \Omega \B \Po \to \Po$ and the Koszul resolution $g_\kappa : \Omega \Poa \to \Po$. These maps do not commute with the semi-augmentations in general.

For example considering the case $\Po = \uAs$ and the associated morphism $g_\kappa : \Omega \uAs^{\ash} \to \uAs$, we can compute
\[
g_\kappa\left(\munsplit - |\right) = 0
\]
whereas $\varepsilon (\munsplit - |) = -|$ and $\varepsilon'(0)=0$.

The consequence is that we have to modify the morphisms between curved coproperads for Theorem 3.4.1 to be correct. 

Unless otherwise stated, we use the notations appearing in the original article.

\section{New conventions}

\subsection{Curvature equation}

First we slightly modify the definition of curved coproperad given in the original article. The only difference is that the sign is changed in the right-hand side of Condition (a) in the following definition.

\begin{dei}
A \emph{curved coproperad} is a triple $(\C,\, d_{\C},\, \theta)$, where $\C$ is a graded coproperad, the map $d_{\C}$ is a coderivation of $\C$ of degree $-1$ and the \emph{curvature} $\theta : \C \rightarrow I$ is a map of degree $-2$ such that:
\begin{enumerate}[label=(\alph*)]
\item ${d_{\C}}^{2} = (\theta \otimes id_{\C} - id_{\C} \otimes \theta) \cdot \Delta_{(1,1)}$,
\item $\theta \cdot d_{\C} = 0$.
\end{enumerate}
\end{dei}

It follows that in Lemma 3.2.2, the same sign has to be changed accordingly without any other changes.

Similarly the same modification has to be done in the definition of curved Lie algebra.

\begin{dei}
A \emph{curved Lie algebra} is a quadruple $(\mathfrak{g},\, [-,\, -],\, d_{\mathfrak{g}},\, \theta)$, where $(\mathfrak{g},\, [-,\, -])$ is a Lie algebra, the map $d_{\mathfrak{g}}$ is a derivation of $\mathfrak{g}$ of degree $-1$ and the curvature $\theta$ is an element of $\mathfrak{g}$ of degree $-2$ such that:
\begin{enumerate}[label=(\alph*)]
\item $d_{\mathfrak{g}}^{2} = [- ,\, \theta]$;
\item $d_{\mathfrak{g}} (\theta) = 0$.
\end{enumerate}
\end{dei}

The sign in Proposition 3.2.4 is changed accordingly.

\subsection{Morphism}

We now modify the definitions of morphism of curved coproperads and of sdg properads as follows.

\begin{dei}
\begin{itemize}
\item
A \emph{lax morphism between curved coproperads} $(\C,\, d_{\C}, \theta) \rightarrow (\C',\, d_{\C'}, \theta')$ is a pair $(f, a)$ where $f : \C \rightarrow \C'$ is a morphism of coproperads and $a : \C \to I$ is an $\So$-bimodule map of degree $-1$ such that
\begin{align}
d_{\C'} \cdot f & = f \cdot d_{\C} + (f \otimes a - a \otimes f) \cdot \Delta_{(1,1)} \quad \textrm{ and} \label{eq: lax morphism diff}\\
\theta' \cdot f & = \theta + a \cdot d_\C + a\star a, \label{eq: lax morphism curvature}
\end{align}
with $a\star a := \gamma_I \cdot(a \otimes a) \cdot \Delta_{(1, 1)}$. 

The composition of lax morphisms is given by $(g, b) \cdot (f, a) := (g\cdot f, a+b\cdot f)$ and the identity is $(\id, 0)$. 
We denote this category by \textsf{curved coprop.} and its morphisms by $\Hom^{\textrm{lax}}$.
\item
A \emph{morphism between two sdg properads} $(\Po,\, d_{\Po},\, \varepsilon) \xrightarrow{f} (\Po',\, d_{\Po'},\, \varepsilon')$ is simply a morphism of dg properads $f : (\Po,\, d_{\Po}) \rightarrow (\Po',\, d_{\Po'})$. 
We denote by \textsf{sdg prop.} the category of semi-augmented dg properads.
\end{itemize}
\end{dei}

\begin{rem}
Morphisms of curved coproperads defined in \cite{HirshMilles} are examples of lax morphisms of curved coproperads with the convention that a morphism $f$ corresponds to the lax morphism $(f, 0)$.
\end{rem}

\section{Corrected results}

\subsection{Curved twisting morphisms}

We prove that for a curved coproperad $\C$ and a dg properad $\Po$, the set of curved twisting morphisms Tw$(\C,\, \Po)$ in Hom$_{\So}(\C,\, \Po)$ forms a bifunctor.\\

Given a coaugmented curved coproperad $(\C,\, d_{\C},\, \theta)$ and a dg properad $(\Po,\, d_{\Po})$, an element $\alpha : \C \rightarrow \Po$ of degree $-1$ in the curved Lie algebra Hom$_{\So}(\C ,\, \Po)$ is called a \emph{curved twisting morphism} if it satisfies that the composition $I \mono \C \xrightarrow{\alpha} \Po$ is equal to $0$ and it is a solution of the \emph{curved Maurer-Cartan equation}
$$\partial(\alpha) + \alpha \star \alpha = \varTheta.$$

\begin{lem}
Curved twisting morphisms forms a bifunctor
\[
\mathrm{Tw}(-,-) : \mathsf{curved\ coprop.}^{\mathrm{op}} \times \mathsf{dg\ prop.} \to \mathsf{Set}.
\]
\end{lem}

\begin{pf}
The functoriality in the right variable is given by post-composition and doesn't cause any trouble. 
We only prove the functoriality in the left variable. Let $(f, a) : \C \to \C'$ be a lax morphism between curved coproperads. To a curved twisting morphism $\alpha : \C' \to \Po$, we associate the map $F(f, a)(\alpha) := \alpha \cdot f + e\cdot a : \C \to \Po$. 
We first compute
\[
\partial (\alpha \cdot f + e\cdot a) = (\partial \alpha)\cdot f - \alpha \cdot (f \otimes a - a \otimes f) \cdot \Delta_{(1, 1)} + e\cdot a \cdot d_\C.
\]
Then
\begin{align*}
(\alpha \cdot f + e\cdot a) \star (\alpha \cdot f + e\cdot a) = (\alpha \star \alpha)\cdot f + \alpha \cdot (f \otimes a - a \otimes f) \cdot \Delta_{(1, 1)} + e\cdot (a \otimes a)\cdot \Delta_{(1, 1)}
\end{align*}
It follows that
\[
\partial (F(f, a)) + F(f, a) \star F(f, a) = e \cdot \theta' \cdot f + e\cdot \left( a \cdot d_\C + (a \otimes a)\cdot \Delta_{(1, 1)}\right) = \Theta.
\]
Therefore $\alpha \cdot f + e\cdot a$ is a curved twisting morphism. 
It is direct to check that $F(\id, 0) = \id$ and $F$ preserves the composition. This concludes the proof.
\end{pf}

\subsection{Bar construction}

We now work in the context of semi-augmented properads. 
The definition of the bar construction is unchanged. 
We correct the proof of Lemma 3.3.4 in \cite{HirshMilles}.

\begin{lem}
The bar construction is a functor $\B : \mathsf{sdg\ prop.} \rightarrow \mathsf{conil.\ curved\ coprop.}$.
\end{lem}

\begin{pf}
Let $f : (\Po,\, d_{\Po},\, \varepsilon) \rightarrow (\Po',\, d_{\Po'},\, \varepsilon')$ be a morphism of sdg properads. The map of dg $\So$-bimodules underlying $f$ is characterized by the morphism of dg $\So$-bimodules $\bar{f} : \overline{\Po} \rightarrow \overline{\Po'}$ and the morphism $\varepsilon' \cdot (s^{-1}f) : s\overline{\Po} \to I$. The map $\F^{c}(\bar f) : \F^{c}(s\overline{\Po}) \rightarrow \F^{c}(s\overline{\Po'})$ is a map of coproperads by construction and $\varepsilon' \cdot (s^{-1}f)$ induces a map $a : \F^{c}(s\overline{\Po}) \to I$ by precomposition by the projection to $s\overline{\Po}$. 
We now show that the couple $(\F^{c}(\bar f), a)$ is a lax morphism of curved coproperads. 

First, the morphism $\bar f$ commutes with $d_{\overline{\Po}}$ and $d_{\overline{\Po'}}$, thus $\F^{c}(\bar f)$ commutes with the coderivations $d_1$ and $d_1'$. 
Then a computation gives the equality on $\F^c(s\overline{\Po})^{(2)}$
\begin{equation}
\label{eq: differentials and lax morphisms}
d_2' \cdot \F^c(\bar f) - \F^c(\bar f) \cdot d_2 = (\bar f\otimes a - a\otimes \bar f)\cdot \Delta_{(1, 1)}.
\end{equation}
Having this equality on $\F^c(s\overline{\Po})^{(2)}$ is equivalent to having the corestriction of this equality to $s\overline{\Po}$. 
Remark moreover that the coproperad $\F^c(s\overline{\Po})$ can be seen as an infinitesimal $\F^c(s\overline{\Po'})$-bicomodule by means of the (coaugmented) coproperad morphism $\F^c(\bar f)$. 
Then a direct computation shows that the map $d_2' \cdot \F^c(\bar f) - \F^c(\bar f)\cdot d_2$ is a coderivation $\F^c(s\overline{\Po}) \to \F^c(s\overline{\Po'})$. (This is a general result for two coderivations $d_2$ and $d_2'$ and a cooperad morphism $\F^c(\bar f)$.)
Similarly, the map $(\F^c(\bar f) \otimes a - a \otimes \F^c(\bar f)) \cdot \Delta_{(1, 1)}$ is a coderivation $\F^c(s\overline{\Po}) \to \F^c(s\overline{\Po'})$ by means of a similar computation as the one given in the proof of Lemma 3.2.2 in \cite{HirshMilles}. 
It follows that Equality \eqref{eq: differentials and lax morphisms} is true on the whole $\F^c(s\overline{\Po})$ (by means of Lemma 15 in \cite{MerkulovVallette} which guarantees that a coderivation is characterized by its corestriction to the generators of the coproperad in the codomain) and that the predifferentials $d_{bar}$ and $d_{bar}'$ satisfy the first equation that a lax morphism of curved coproperads should satisfy (that is \eqref{eq: lax morphism diff}). 

For the second equation (that is \eqref{eq: lax morphism curvature}), we can easily check that the equality $\theta'_{bar} \cdot \F^{c}(\bar f) = \theta_{bar}$ is satisfied on $s\overline{\Po}$. 
On $\F^c(s\overline{\Po})^{(2)}$, we have
\[
\theta'_{bar} \cdot \F^{c}(\bar f) - \theta_{bar} = a\cdot d_{bar} + \gamma_I \cdot(a\otimes a)\cdot \Delta_{(1, 1)}.
\]
Finally the second equation is satisfied and we obtain that the bar construction is a functor.
$\cqfd$
\end{pf}

\subsection{Cobar construction}

We work in the context of coaugmented curved coproperads. 
The definition of the cobar construction is unchanged. 
We correct the proof of Lemma 3.3.7 in \cite{HirshMilles}.

\begin{lem}
The cobar construction is a functor $\Omega : \mathsf{coaug.\ curved\ coprop.} \rightarrow \mathsf{sdg\ prop.}$.
\end{lem}

\begin{pf}
Let $(f, a) : (\C,\, d_{\C},\, \theta) \rightarrow (\C',\, d_{\C'},\, \theta')$ be a lax morphism between two coaugmented curved coproperads. We recall that $f$ has degree 0 (as a morphism of coproperads). 
We define the map $\F(f, a) : \F(s^{-1}\overline{\C}) \rightarrow \F(s^{-1}\overline{\C'})$ by its restriction to $s^{-1}\overline{\C}$ as the following formula
\[
s^{-1}\overline{\C} \xrightarrow{-sa  + f} I \oplus s^{-1}\overline{\C} \ \subset\ \F(s^{-1}\overline{\C}).
\]
(Because $f$ is a morphism of coaugmented coproperad, it sends $\overline{\C}$ to $\overline{\C'}$.) 
Noting $d$, resp. $d'$, the differential on $\F(s^{-1}\overline{\C})$, resp. on $\F(s^{-1}\overline{\C'})$, the map $d' \cdot \F(f, a) - \F(f, a) \cdot d$ is a derivation. To prove that it is zero, it is therefore enough to prove that it is zero on the generators $s^{-1}\overline{\C}$ (\cite[Lemma 14]{MerkulovVallette}). 
We have on $s^{-1}\overline{\C}$
\begin{align}
-d_0' \F(f, a) +  \F(f, a) d_0 & = -d_0' \cdot (-sa + f) + \F(f, a) \cdot (s\theta) = -s(\theta' \cdot f -\theta)\\
d_1' \F(f, a) -  \F(f, a) d_1 & = d_1'(-sa+f)-\F(f, a) d_{s^{-1}\overline{\C}} = \id_{s^{-1}}\otimes (d_{\overline{\C'}} f - f d_{\overline{\C}}) + s ad_{s^{-1}\overline{\C}}\\
-d_2' \F(f, a) +  \F(f, a) d_2 & = -d_2'(-sa+f)+\F(f, a) (s^{-1}\overline{\Delta}_{(1,1)})\\
& = -s^{-1}\overline{\Delta'}_{(1,1)} f + (-sa+f) \otimes (-sa+f) s^{-1}\overline{\Delta}_{(1, 1)} \notag\\
& = \left(s^{-1}(f \otimes a - a\otimes f) + a\otimes a\right)\Delta_{(1, 1)}(s -). \notag
\end{align}
It follows from these three equations, using the definition of lax morphisms of curved coproperads, that $d' \cdot \F(f) = \F(f) \cdot d$.
$\cqfd$
\end{pf}

\subsection{Bar-cobar adjunction}

We finally prove that the bar and the cobar functors form an adjoint pair.

\begin{thm}\label{barcobaradjunction}
For any conilpotent curved coproperad $\C$ and for any sdg properad $\Po$,
there are natural bijections
$$\mathrm{Hom}_{\mathsf{sdg\ prop.}}(\Omega \C,\, \Po) \cong \mathrm{Tw}(\C,\, \Po) \cong \mathrm{Hom}_{\mathsf{coaug.\ curved\ coprop.}}^{\mathrm{lax}}(\C,\, \B \Po).$$
\end{thm}

\begin{pf}
We make the first bijection explicit. A morphism of semi-augmented (graded) properads $f_{\alpha} : \F(s^{-1}\overline{\C}) \rightarrow \Po$ is uniquely determined by a map $s \otimes \alpha : s^{-1}\overline{\C} \rightarrow \Po$ of degree $0$, or equivalently, by a map $\alpha : \C \rightarrow \Po$ of degree $-1$ satisfying $I \mono \C \xrightarrow{\alpha} \Po$ is zero.

Moreover, $f_{\alpha}$ commutes with the differentials if and only if the following diagram commutes
$$\xymatrix{s^{-1}\overline{\C} \ar[r]^{s \otimes \alpha} \ar[d]_{-d_{0} + d_{1} - d_{2}} & \Po \ar[r]^{d_{\Po}} & \Po\\
\F(s^{-1}\overline{\C}) \ar[rr]_{\F(s \otimes \alpha)} && \F(\Po), \ar[u]_{\widetilde{\gamma}}}$$
where $\widetilde{\gamma}$ is induced by $\gamma$. We have
$$\begin{array}{l}
d_{\Po} \cdot (s \otimes \alpha) = - s \otimes (d_{\Po} \cdot \alpha),\\
\widetilde{\gamma} \cdot \F(s \otimes \alpha) \cdot d_0 = e \cdot (s \otimes \theta) = s \otimes (e \cdot \theta),\\
\widetilde{\gamma} \cdot \F(s \otimes \alpha) \cdot d_1 = (s \otimes \alpha) \cdot (id_{s^{-1}} \otimes d_{\C}) = s \otimes (\alpha \cdot d_{\C}),\\
\widetilde{\gamma} \cdot \F(s \otimes \alpha) \cdot d_2 = \gamma \cdot \big((s \otimes \alpha) \otimes (s \otimes \alpha)\big) \cdot (s^{-1} \otimes \overline{\Delta}_{(1,\, 1)}) = - s \otimes (\gamma \cdot (\alpha \otimes \alpha) \cdot \Delta_{(1,\, 1)}).
\end{array}$$
Thus the commutativity of the previous diagram is equivalent to the equality
$$-e \cdot \theta + \alpha \cdot d_{\C} + \gamma \cdot (\alpha \otimes \alpha) \cdot \Delta_{(1,\, 1)} = -d_{\Po} \cdot \alpha,$$
that is $\partial(\alpha) + \alpha \star \alpha = \varTheta$.

We now make the second bijection explicit. A morphism of coaugmented coproperads $g_{\alpha} : \C \rightarrow \F^{c}(s\overline{\Po})$ is uniquely determined by a map $s \otimes \bar\alpha : \C \rightarrow s\overline{\Po}$ which sends $I$ to $0$, that is by a map $\bar\alpha : \C \rightarrow \Po$ of degree $-1$ satisfying $I \mono \C \xrightarrow{\bar\alpha} \Po$ and $\C \xrightarrow{\bar\alpha} \Po \xrightarrow{\varepsilon} I$ are zero. 
It follows that a lax morphism between curved coproperads $(g_\alpha, a_\alpha) : \C \to \B\Po$ can be written in a unique way as a map $\alpha = \bar\alpha + e\cdot a_\alpha : \C \rightarrow \Po$ of degree $-1$ satisfying $I \mono \C \xrightarrow{\alpha} \Po$ is zero (by using that $a_{\alpha}$ has degree $-1$).

Moreover, $(g_{\alpha}, a_\alpha)$ satisfies Equations \eqref{eq: lax morphism diff} and \eqref{eq: lax morphism curvature} if and only if the following diagrams commute up to the terms $(g_\alpha \otimes a_\alpha - a_\alpha \otimes g_\alpha) \cdot \Delta_{(1, 1)}$ and $a_\alpha \cdot d_\C + (a_\alpha \otimes a_\alpha)\cdot \Delta_{(1,1)}$ respectively
$$\xymatrix@C=120pt@R=10pt{\C \ar[r]^{s \otimes \bar\alpha + \big((s \otimes \bar\alpha) \otimes (s \otimes \bar\alpha)\big) \cdot \Delta_{(1,\, 1)} \hspace{1cm}} \ar[dd]_{d_{\C}} & s\overline{\Po} \oplus s\overline{\Po} \boxtimes_{(1,\, 1)} s\overline{\Po} \ar[dd]^{d_{bar} = d_{1} + d_{2}}\\
&\\
\C \ar[r]_{s \otimes \bar\alpha} & s\overline{\Po}}
\hspace{1.5cm} \textrm{and} \hspace{1.5cm} \mbox{\xymatrix@R=10pt{\C \ar[dd]_{\theta} \ar[r]^{g_{\alpha}} & \B \Po \ar[ddl]^{\theta_{bar}}\\
&\\
I. &}}$$
Since $\bar\alpha \star \bar\alpha = -(s^{-1} \otimes inc) \cdot d_{2} \cdot \big((s \otimes \bar\alpha) \otimes (s \otimes \bar\alpha)\big) \cdot \Delta_{(1,\, 1)} + e \cdot \theta_{bar} \cdot g_{\alpha}$, the commutativity of the diagrams (up to the above corresponding terms) gives $\partial(\alpha) + \alpha \star \alpha = \varTheta$. Moreover, the projections of the curved Maurer-Cartan equation on $\overline{\Po}$ and on $I$ give the two commutative diagrams (up to the above corresponding terms). This concludes the proof.
$\cqfd$
\end{pf}

\begin{ex}
\begin{itemize}
\item[]
\item To the identity morphism $(id_{\B \Po}, 0) : \B \Po \rightarrow \B \Po$ of coaugmented curved coproperads corresponds the curved twisting morphism $\pi : \B \Po \rightarrow \Po$ defined by $\F^{c}(s\overline{\Po}) \epi s\overline{\Po} \cong \overline{\Po} \mono \Po$.
\item To the identity morphism $id_{\Omega \C} : \Omega \C \rightarrow \Omega \C$ of properads corresponds the curved twisting morphism $\iota : \C \rightarrow \Omega \C$ defined by $\C \rightarrow \overline{\C} \cong s^{-1}\overline{\C} \mono \F(s^{-1}\overline{\C})$.
\end{itemize}
\end{ex}

\begin{lem}\label{universaltwmorph}
For any conilpotent curved coproperad $\C$ and for any sdg properad $\Po$, every curved twisting morphism $\alpha : \C \rightarrow \Po$ factors through the universal curved twisting morphisms $\pi$ and $\iota$:
$$\xymatrix{& \Omega \C \ar@{-->}[dr]^{-f_{\alpha}} &\\
\C \ar[rr]^{\alpha} \ar[ur]^{\iota} \ar@{-->}[dr]_{(g_{\alpha}, a_\alpha)} && \Po\\
& \B \Po, \ar[ur]_{\pi} &}$$
where $f_{\alpha}$ is a morphism of sdg properads and $(g_{\alpha}, a_\alpha)$ is a lax morphism of conilpotent curved coproperads (with $a_\alpha := \varepsilon \cdot \alpha$). 
The bottom factorization means that the diagram is commutative up to $e\cdot a_\alpha$.
\end{lem}

\begin{pf}
The dashed arrows are just the images of $\alpha$ under the two bijections of Proposition \ref{barcobaradjunction}. The diagram commutes by a direct calculation.
$\cqfd$
\end{pf}

The other results of \cite{HirshMilles}, with the convention that a morphism $f$ corresponds to the lax morphism $(f, 0)$, are unchanged.

\bibliographystyle{plain}
\bibliography{bib}
\ \\
\noindent
\textsc{Joan Mill\`es, Université Paul Sabatier, Institut de mathématiques de Toulouse, 118, route de Narbonne, F-31062 Toulouse Cedex 9, France}\\
E-mail : \texttt{joan.milles@math.univ-toulouse.fr}\\

\noindent
\textsc{Joseph Hirsh}\\
E-mail : \texttt{josephhirsh@gmail.com}
\end{document}